\documentclass[10pt,a4paper]{article}
\title{\textbf{Mixed $3$-Sasakian structures and curvature}}
\author{Angelo V. Caldarella and Anna Maria Pastore}
\date{ }

\usepackage{amsfonts}
\usepackage{amsmath,amssymb,amscd,amsthm}
\usepackage{multicol,longtable,tabularx}
\usepackage{enumerate}
\usepackage{latexsym}
\usepackage{amstext,graphics,color}
\usepackage{ifthen}
\usepackage{verbatim}
\usepackage{calc}

\newtheorem{theorem}{Theorem}[section]

\newtheorem{lemma}[theorem]{Lemma}
\newtheorem{proposition}[theorem]{Proposition}
\theoremstyle{definition}

\newtheorem{definition}[theorem]{Definition}
\newtheorem{example}[theorem]{Example}

\newtheorem{remark}[theorem]{Remark}

\newcommand{\proofend}{\hfill{\mbox{$\Box$}}}
\newcommand{\proofbegin}{\noindent{\sc Proof.\ }}
\renewcommand{\proof}{\proofbegin}
\renewcommand{\endproof}{\proofend}
\newcommand{\pmm}{$(M,\varphi,\xi,\eta,g)$}
\newcommand{\mts}{$(\varphi_a,\xi_a,\eta^a)$}
\newcommand{\mmts}{$(\varphi_a,\xi_a,\eta^a,g)$}

\newcommand{\R}{\mathbb{R}}

\def\stackunder#1#2{\displaystyle{\mathrel{\mathop{#2}\limits_{#1}}}}

\begin{document}
\maketitle
\begin{abstract}
In this paper we deal with two classes of mixed metric 3-structures, namely the mixed 3-Sasakian structures and the mixed metric 3-contact structures. Firstly we study some properties of the curvature of mixed 3-Sasakian structures, proving that any manifold endowed with such a structure is Einstein. Then we prove the identity between the class of mixed 3-Sasakian structures and the class of mixed metric 3-contact structures.
\end{abstract}
\textbf{2000 Mathematics Subject Classification} 53C15, 53C25, 53C50.\\
\textbf{Keywords and phrases.} Mixed metric 3-structures; mixed 3-Sasakian structures; mixed metric 3-contact structures; paracontact structures; semi-Riemannian manifolds; Einstein manifolds.

\section{Introduction}
The geometry of $3$--Sasakian manifolds is a well-known topic, since its introduction, made indipendently by Kuo (\cite{Ku}) and Udri\c{s}te (\cite{U}), which has been extensively and deeply studied, in a first stage, by Ishihara, Kashiwada, Konishi, Kuo, Tachibana, Tanno, Yu and other geometers of the Japanese school, and then, from a different viewpoint, by Boyer, Galicki and Mann. Among the various works and papers, we only recall the last chapter of \cite{Bla} and the remarkable survey \cite {BG}, to which we refer the reader for more details about such structures, as well as for historical remarks and bibliographic references. On the other hand, only very recently a first study of the analogue odd-dimensional geometries related to the algebra of paraquaternionic numbers has been conducted by Ianu\c s, Mazzocco and V\^\i lcu in \cite{IMV} and \cite{IV}.

In analogy with an early result of Kashiwada (\cite{Ka01}) for Sasakian 3-structures, the first of the two main results we shall present in this paper is for manifolds endowed with mixed 3-Sasakian structures. Namely, any manifold with a mixed 3-Sasakian structure is Einstein. This can be also regarded as an analogous to the well-known fact that a paraquaternionic K\"ahler manifold is Einstein (cf.\ \cite{G-RMV-L}). To this purpose, we shall need some formulas for the curvature tensor of a manifold with paraSasakian structure and of a manifold with indefinite Sasakian structure. Some formulas recently proved in \cite{Z} shall be here recovered, even if in a slightly more general way, because of the fact that we consider here both the cases of positive and negative paraSasakian structures.

The second main result is concerned with the identity between the class of mixed metric 3-contact structures and the class of mixed 3-Sasakian structures (see Kashiwada \cite{Ka02} for the case of 3-contact metric manifolds). Its achievement is based on an extension of Kashiwada's generalization of a lemma of Hitchin (cf. \cite{Ka03}) to the almost hyper paraHermitian case.

The content of the paper is now briefly described.

In section $2$ we give some fundamental definitions and facts about paracontact metric structures (cf.\ \cite{E02}, \cite{Z}), which, together with the notion of indefinite almost contact metric structure (\cite{B}), are at the root of the notion of mixed metric 3-structure. We also recall few definitions concerning almost hyper paraHermitian structures. In section $3$, after introducing the notation of $[r]$-Sasakian structure, $r=\pm1$, to mean an indefinite Sasakian structure for $r=+1$, and a paraSasakian structure for $r=-1$, we prove some preliminary issues, needed to state, in section $4$, the result concerning the mixed 3-Sasakian manifolds. Finally, section $5$ is devoted to prove that a mixed metric 3-contact structure is in fact a mixed 3-Sasakian structure.

All manifolds and tensor fields are assumed to be smooth.\\

\section{Preliminaries}

We recall few definitions about paracomplex and hyper paracomplex structures. For more details we refer to \cite{CFG} and \cite{IZ}.

\begin{definition}
An \emph{almost product structure} on a manifold $M$ is a $(1,1)$--type tensor field $F\neq\pm I$, satisfying $F^{2}=I$; the pair $(M,F)$ is then said to be an \emph{almost product manifold}.
\end{definition}

On an almost product manifold $(M,F)$ we have $TM=T^+M\oplus T^-M$, where $T^+M$ and $T^-M$ are the eigensubbundles related to the eigenvalues $+1$ and $-1$ of $F$. If $\mathrm{rank}(T^+M)=\mathrm{rank}(T^-M)$, then $(M,F)$ is said to be \emph{almost paracomplex manifold}. Finally, an almost product (paracomplex) manifold $(M,F)$ is called \emph{product (paracomplex) manifold}, if $N_F=0$, being $N_F$ the Nijenhuis tensor field of the structure $F$. Any (almost) paracomplex manifold has even dimension.

An (almost) paracomplex manifold $(M,F)$ is called \emph{(almost) paraHermitian manifold} if there exists a metric tensor $g$ compatible with $F$, i.e. such that $g(FX,Y)+g(X,FY)=0$, for any $X,Y\in\Gamma(TM)$. Such a metric is necessarily semi-Riemannian, with neutral signature.

\begin{definition}\label{def005}
An \emph{almost hyper paraHermitian structure} on a manifold $M
$ is a triple $(J_{1},J_{2},J_{3})$ of $(1,1)$--type tensor fields, together with a semi--Riemannian metric $g$
satisfying:
\begin{enumerate}
	\item[(i)] $(J_a)^2=-\tau_aI$, \quad for any $a\in\{1,2,3\}$,
	\item[(ii)] $J_aJ_b=\tau_cJ_c=-J_bJ_a$, \quad for any cyclic permutation $(a,b,c)$ of $(1,2,3)$,
	\item[(iii)] $g(J_aX,Y)+g(X,J_aY)=0$, \quad for any $a\in\{1,2,3\}$, and any $X,Y\in\Gamma(TM)$,
\end{enumerate}
where $\tau_1=-1$, $\tau_2=-1$ and $\tau_3=+1$. Then $(M,J_1,J_2,J_3,g)$ will be said to be an \emph{almost hyper paraHermitian manifold}.
\end{definition}

Such a manifold has dimension divisible by four, and the metric has neutral signature. An almost hyper paraHermitian structure on a manifold $M$ will be said \emph{hyper paraHermitian} if and only if, for any $a\in\{1,2,3\}$, the Nijenhuis tensor field $N_a$ vanishes, that is each structure $J_a$ is integrable. Then $M$ will be called \emph{hyper paraHermitian manifold}. An almost hyper paraHermitian manifold is hyper paraHermitian if and only if at least two of the Nijenhuis tensor fields vanish (cf. \cite{IZ}).

\begin{definition}
Let $M$ be a manifold. An \emph{almost paracontact structure} on $M$ is a system $(\varphi,\xi,\eta)$, where $\varphi\in\mathfrak{T}^1_1(M)$, $\xi\in\Gamma(TM)$ and $\eta\in\bigwedge^1(M)$, satisfying $\varphi^2=I-\eta\otimes\xi$ and $\eta(\xi)=1$. Then $M$ is said to be an \emph{almost paracontact manifold}, denoted by $(M,\varphi,\xi,\eta)$. An almost paracontact structure $(\varphi,\xi,\eta)$ will be said \emph{normal} if $N_\varphi=2d\eta\otimes\xi$, where $N_\varphi$ is the Nijenhuis tensor field of the structure $\varphi$.
\end{definition}

Almost paracontact structures were originally introduced by I. Sat\=o in \cite{Sa01} and \cite{Sa02}, where he studied also the properties of manifolds endowed with such structures and with a Riemannian metric satisfying suitable compatibility conditions. Besides, one may find similar definitions in \cite{KW} and \cite{Z}, where the further condition that the restriction $\varphi|_{\mathrm{Im}(\varphi)}$ is an almost paracomplex structure on the distribution $\mathrm{Im}(\varphi)$ is required. The notion of normality for an almost paracontact structure is defined, as in the classical almost contact case (cf.\ \cite{Bla}), through the integrability of the almost product structure $F$ canonically induced on the manifold $M\times\R$, defined by $F\left(X,f\frac{d}{dt}\right):=\left(\varphi X+f\xi,\eta(X)\frac{d}{dt}\right)$ (cf. \cite{KW}, \cite{Z}).

Other properties of an almost paracontact manifold $(M,\varphi,\xi,\eta)$, which are immediate consequences of the above definition, are $\varphi(\xi)=0$, $\eta\circ\varphi=0$, $\ker(\varphi)=\mathrm{Span}(\xi)$, $\ker({\eta})=\mathrm{Im}(\varphi)$ and $TM=\mathrm{Im}(\varphi)\oplus\mathrm{Span}(\xi)$.

Endowing an almost paracontact manifold with a metric tensor field, and considering a suitable compatibility condition, we obtain the notion of almost paracontact metric manifold.

\begin{definition}[\cite{Z}]
Let $(M,\varphi,\xi,\eta)$ be an almost paracontact manifold, and $g$ a metric tensor field on $M$, that is a symmetric, non degenerate $(0,2)$--type tensor field on $M$. Then $g$ is said to be \emph{compatible} with the structure $(\varphi,\xi,\eta)$, if the following condition holds
\begin{equation*}
g(\varphi X,\varphi Y)=-g(X,Y)+\varepsilon\eta(X)\eta(Y)
\end{equation*}
for any $X,Y\in\Gamma(TM)$, with $\varepsilon=\pm1$, according as $\xi$ is spacelike or timelike, respectively. Then, the structure $(\varphi,\xi,\eta,g)$ is said to be an \emph{almost paracontact metric structure}. We shall call the structure \emph{positive} or \emph{negative}, according as $\varepsilon=+1$ or $\varepsilon=-1$, respectively. Then $(M,\varphi,\xi,\eta,g)$ will be called \emph{almost paracontact metric manifold}. Such a structure $(\varphi,\xi,\eta,g)$ will be called \emph{normal} if $N_\varphi=2d\eta\otimes\xi$.
\end{definition}

In \cite{E02}, the author refers to the same kind of structure with the denomination of almost paracontact hyperbolic metric structure.

As a consequence of the above definition, for an almost paracontact metric manifold \pmm{ }one has that the pair $(F,g)$, where $F:=\varphi|_{\mathrm{Im}(\varphi)}$, is an almost paraHermitian structure on the distribution $\mathrm{Im}(\varphi)$, hence $\mathrm{rank}({\mathrm{Im}(\varphi)})=2m$ and $\dim(M)=2m+1$. Furthermore, the signature of $g$ on $\mathrm{Im}(\varphi)$ is $(m,m)$, where we put first the minus signs, and the signature of $g$ on $TM$ is $(m,m+1)$ or $(m+1,m)$, according as $\xi$ is spacelike (the structure is positive) or timelike (the structure is negative), respectively. It follows that $g$ is a Lorentzian metric only if $m=1$ and $\dim(M)=3$.

We have that $TM$ is the orthogonal direct sum of $\mathrm{Im}(\varphi)$ and $\mathrm{Span}(\xi)$, and finally that $\eta(X)=\varepsilon g(X,\xi)$ and $g(\varphi X,Y)+g(X,\varphi Y)=0$, for any $X,Y\in\Gamma(TM)$.

Particular classes of almost paracontact metric structures are defined as follows.

\begin{definition}[\cite{E02}, \cite{Z}]
Let $(M,\varphi,\xi,\eta,g)$ be an almost paracontact metric manifold. Then it is said to be a
\begin{enumerate}
  \item [(i)] \emph{paracontact metric manifold}, if $d\eta=\Phi$;
	\item [(ii)] \emph{paraSasakian manifold}, if $d\eta=\Phi$ and the structure is normal;
	\item [(iii)] \emph{para--K--contact manifold}, if $d\eta=\Phi$ and $\xi$ is a Killing vector field,
\end{enumerate}
being $\Phi(X,Y):=g(X,\varphi Y)$ the fundamental 2-form associated with the almost paracontact metric structure.
\end{definition}

Furthermore, we recall the following result.
\begin{proposition}\label{prop021}
Let $(M,\varphi,\xi,\eta,g)$ be an almost paracontact metric manifold. Then, it is a paraSasakian manifold if and only if, setting $\varepsilon=g(\xi,\xi)=\pm1$, 
\begin{equation*}
(\nabla_X\varphi)(Y)=-g(X,Y)\xi+\varepsilon\eta(Y)X
\end{equation*}
for any $X,Y\in\Gamma(TM)$.
\end{proposition}
 
We assume the following definition of mixed (metric) 3--structure, which, althought in a different form, is introduced in \cite{IMV} and \cite{IV}.

\begin{definition}
Let $M$ be a manifold. A \emph{mixed 3-structure} on $M$ is a triple of structures \mts,\ $a\in\{1,2,3\}$, which are almost paracontact structures for $a=1,2$ and almost contact structure for $a=3$, and satisfy
{\setlength\arraycolsep{0pt}
\begin{eqnarray}
&&\varphi_{a}\varphi_{b}-\tau_{a}\eta^{b}\otimes\xi_{a}=\tau_{c}\varphi_{c}=-\varphi_{b}\varphi_{a}+\tau_{b}\eta^{a}\otimes\xi_{b}, \label{eq90}\\
&&\eta^{a}\circ\varphi_{b}=\tau_{c}\eta^{c}=-\eta^{b}\circ\varphi_{a}, \label{eq91}\\
&&\varphi_{a}(\xi_{b})=\tau_{b}\xi_{c}, \quad\quad \varphi_{b}(\xi_{a})=-\tau_{a}\xi_{c}, \label{eq92}
\end{eqnarray}}for any cyclic permutation $(a,b,c)$ of $(1,2,3)$, with $\tau_1=\tau_2=-1=-\tau_3$. If, further, there exists on $M$ a metric tensor field $g$ satisfying
\begin{equation}\label{eq93}
g(\varphi_{a}X,\varphi_{a}Y)=\tau_{a}\left(g(X,Y)-\varepsilon_{a}\eta^{a}(X)\eta^{a}(Y) \right),
\end{equation}
for any $a\in\{1,2,3\}$, and any $X,Y\in\Gamma(TM)$, where $\varepsilon_a=g(\xi_a,\xi_a)=\pm1$, then the metric tensor $g$ is said to be \emph{compatible} with the mixed 3-structure, and \mmts,\ $a\in\{1,2,3\}$, is called a \emph{mixed metric 3-structure}.
\end{definition}

From now on, a mixed 3-structure and a mixed metric 3-structure on a manifold $M$ will be denoted simply with \mts{ }and \mmts,{ }respectively, leaving the condition $a\in\{1,2,3\}$ understood.

\begin{remark}
We point out that the above definition of mixed 3-structure, without the metric compatibility, is equivalent to the definition given in \cite{IMV}, and very recently in \cite{IV}, providing that one substitutes the structures $(\varphi_1,\xi_1,\eta^1)$, $(\varphi_2,\xi_2,\eta^2)$ and $(\varphi_3,\xi_3,\eta^3)$ of \cite{IMV} and \cite{IV} with $(\varphi_3,\xi_3,\eta^3)$, $(\varphi_1,\xi_1,\eta^1)$ and $(\varphi_2,\xi_2,\eta^2)$, respectively, and then the vector fields $\xi_1$ and $\xi_2$ with their opposite. Nevertheless, the metric compatibility expressed in \cite{IV} is also equivalent to (\ref{eq93}) of the above definition only in the case $(\varepsilon_1,\varepsilon_2,\varepsilon_3)=(-1,-1,+1)$.
\end{remark}

Let $M$ be a manifold endowed with a mixed 3-structure \mts. Considering the two distributions $\mathcal{H}:=\bigcap\nolimits_{a=1}^3\ker(\eta^a)$ and $\mathcal{V}:=\mathrm{Span}(\xi_1,\xi_2,\xi_3)$, then one has the decomposition $TM=\mathcal{H}\oplus\mathcal{V}$. It follows that $(\varphi_1|_{\mathcal{H}},\varphi_2|_{\mathcal{H}},\varphi_3|_{\mathcal{H}})$ is an almost hyper paracomplex structure on the distribution $\mathcal{H}$. Hence $\mathrm{rank}(\mathcal{H})=2n$ and $\dim(M)=2n+3$. Furthermore, if we have a mixed metric 3-structure \mmts{ }on $M$, then $(\varphi_a|_{\mathcal{H}},g)$, $a\in\{1,2,3\}$, becomes an almost hyper paraHermitian structure on the distribution $\mathcal{H}$. Hence $\mathrm{rank}(\mathcal{H})=4m$ and $\dim(M)=4m+3$. As an obvious consequence, we have the following result.

\begin{proposition}
Let $M$ be a manifold with $\dim(M)=2n+3$, endowed with a mixed 3-structure \mts. If $n\neq2m$, then there is no metric tensor field $g$ on $M$ compatible with the mixed 3-structure, and $M$ can not have any mixed metric 3-structure.
\end{proposition}

The compatibility condition (\ref{eq93}) between a metric tensor $g$ and a mixed 3-structure \mts{ }on a $(4m+3)$--dimensional manifold $M$, together with (\ref{eq92}), has some consequences also on the signature of the metric $g$. In fact, it is easily checked that $g(\xi_1,\xi_1)=g(\xi_2,\xi_2)=-g(\xi_3,\xi_3)$, hence the vector fields $\xi_1$ and $\xi_2$, related to the almost paracontact metric structures, are both either spacelike or timelike. We may therefore distinguish between \emph{positive} and \emph{negative} mixed metric 3-structures, according as $\xi_1$ and $\xi_2$ are both spacelike ($\varepsilon_1=\varepsilon_2=+1$), or both timelike ($\varepsilon_1=\varepsilon_2=-1$), respectively. This forces the causal character of the third vector field $\xi_3$. Being the signature of $g$ on $\mathcal{H}$ necessarily neutral $(2m,2m)$, we may have only the following two cases
\begin{enumerate}
	\item [(i)] the signature of $g$ on $TM$ is $(2m+1,2m+2)$, if the mixed metric 3-structure is positive $(\varepsilon_1,\varepsilon_2,\varepsilon_3)=(+1,+1,-1)$;
	\item [(ii)] the signature of $g$ on $TM$ is $(2m+2,2m+1)$, if the mixed metric 3-structure is negative $(\varepsilon_1,\varepsilon_2,\varepsilon_3)=(-1,-1,+1)$.
\end{enumerate}
We point out that any metric $g$ which is compatible with a mixed 3-structure, in the sense of (\ref{eq93}), can never be a Lorentzian one and that the definition of mixed metric 3-structure given in \cite{IV} is equivalent to that of a negative mixed metric 3-structure.

\begin{example}[\cite{C}]
Let $M^{4m+3}$ be any orientable non degenerate hypersurface of an almost hyper paraHermitian manifold $(\bar{M}^{4m},J_a,G)_{a=1,2,3}$. If $N\in\Gamma(TM^\bot)$ is a unit normal vector field, such that $G(N,N)=s=\pm1$, one puts $\xi_a:=-\tau_aJ_aN$, for any $a\in\{1,2,3\}$, and defines three $(1,1)$--type tensor fields $\varphi_a$, and three 1--forms $\eta^a$ on $M$, such that $J_aX=\varphi_aX+\eta^a(X)N$, for any $X\in\Gamma(TM)$ and any $a\in\{1,2,3\}$. Then, denoting with $g$ the metric induced on $M$ from $G$, one checks that \mmts{ }is a mixed metric 3--structure on $M$, of sign $\sigma=-s$.
\end{example}

Finally, we introduce the notion of mixed 3-Sasakian structure on a manifold with the following definition.

\begin{definition}\label{Def999}
Let $M$ be a manifold with a mixed metric 3-structure \mmts. This structure will be said a \emph{mixed 3-Sasakian structure} if $(\varphi_1,\xi_1,\eta^1,g)$ and $(\varphi_2,\xi_2,\eta^2,g)$ are paraSasakian structures, and $(\varphi_3,\xi_3,\eta^3,g)$ is an indefinite Sasakian structure. Then $(M,\varphi_a,\xi_a,\eta^a,g)$ will be said \emph{mixed 3-Sasakian manifold}.
\end{definition}

\begin{remark}
By Proposition \ref{prop021}, it follows that a mixed metric 3-structure \mmts{ }on a manifold $M$ is mixed 3-Sasakian if and only if
\begin{equation}\label{eq114}
(\nabla_X\varphi_a)(Y)=\tau_a\left(g(X,Y)\xi_a-\varepsilon_a\eta^a(Y)X\right),
\end{equation}
for any $X,Y\in\Gamma(TM)$ and any $a\in\{1,2,3\}$, with $\tau_1=\tau_2=-1=-\tau_3$.
\end{remark}

We remark that Definition \ref{Def999} is not equivalent to that given in \cite{IV}. More precisely, referring to \cite{IV}, the condition $(\nabla_X\varphi_2)(Y)=g(\varphi_2X,\varphi_2Y)\xi_2+\eta^2(Y)(\varphi_2)^2(X)$ in the Definition 4.3, using the compatibility condition (29), may be rewritten in the form
$(\nabla_X\varphi_2)(Y)=-g(X,Y)\xi_2+\eta^2(Y)X$, which corresponds to 
\begin{equation}\label{eq115}
(\nabla_X\varphi_1)(Y)=g(X,Y)\xi_1+\eta^1(Y)X.
\end{equation}
on our notation.
Since the definition of mixed metric 3-structure given in \cite{IV} is equivalent to that of negative mixed metric 3-structure, writing the condition (\ref{eq114}) for $\tau_a=\tau_1=-1$ and $\varepsilon_a=\varepsilon_1=-1$, we get
$(\nabla_X\varphi_1)(Y)=-g(X,Y)\xi_1-\eta^1(Y)$,
which is clearly the opposite of (\ref{eq115}). One obtains an analogous result considering the condition on $(\nabla_X\varphi_3)(Y)$ of \cite{IV}.  

\section{On the curvature of $[r]$--Sasakian structures}
In this section, we prove some useful formulas concerning the curvature of both paraSasakian structures and indefinite Sasakian structures.
To treat both cases simultaneously, we introduce the synthetic notation of \emph{$[r]$--Sasakian structure} on a manifold $M$, considering a system $(\varphi,\xi,\eta,g)$ where $\varphi\in\mathfrak{T}^1_1(M)$, $\xi\in\Gamma(TM)$, $\eta\in\bigwedge^1(M)$ and $g\in\mathfrak{T}^0_2(M)$ is a metric tensor field, such that $g(\xi,\xi)=\varepsilon=\pm1$, $\varphi^2=r(-I+\eta\otimes\xi)$, $\eta(\xi)=1$ and
{\setlength\arraycolsep{0pt}
\begin{eqnarray}
&&g(\varphi X,\varphi Y)=r(g(X,Y)-\varepsilon\eta(X)\eta(Y)),\label{Eq303} \\
&&(\nabla_X\varphi)(Y)=r(g(X,Y)\xi-\varepsilon\eta(Y)X)\label{eq303},
\end{eqnarray}}obtaining an indefinite Sasakian structure for $r=+1$, and a paraSasakian structure for $r=-1$. From (\ref{eq303}) it follows that $\nabla_X\xi=-\varepsilon\varphi(X)$, for any $X\in\Gamma(TM)$.

Following \cite{KN}, the curvature tensor field $R\in\mathfrak{T}^1_3(M)$ of the Levi-Civita connection $\nabla$, the Riemannian curvature tensor field $R\in\mathfrak{T}^0_4(M)$, and the Ricci curvature tensor field $\rho\in\mathfrak{T}^0_2(M)$ will be defined by
{\setlength\arraycolsep{0pt}
\begin{eqnarray*}
&&R(X,Y)Z:=\nabla_X\nabla_YZ-\nabla_Y\nabla_XZ-\nabla_{[X,Y]}Z,\\
&&R(X,Y,Z,W):=g(R(Z,W)Y,X)=-g(R(X,Y)W,Z),\\
&&\rho(X,Y):=\mathrm{tr}_g\{Z\mapsto R(Z,X)Y\}=\sum_{i=1}^m\varepsilon_ig(R(E_i,X)Y,E_i),
\end{eqnarray*}}where $(E_i)_{1\leqslant i\leqslant m}$ is a local orthonormal frame, $\varepsilon_i=g(E_i,E_i)$ and $m=\dim(M)$.

\begin{lemma}\label{lemma139}
Let $M$ be a manifold endowed with an $[r]$--Sasakian structure $(\varphi,\xi,\eta,g)$. Then, for any $X,Y,Z,W\in\Gamma(TM)$
\begin{equation}\label{eq139}
g(R(X,Y)Z,\varphi W)+g(R(X,Y)\varphi Z,W)=-r\varepsilon P(X,Y,Z,W),
\end{equation}
where $P\in\mathfrak{T}^0_4(M)$ is the tensor field defined by
{\setlength\arraycolsep{2pt}
\begin{eqnarray*}
P(X,Y,Z,W)&:=&d\eta(X,Z)g(Y,W)-d\eta(X,W)g(Y,Z)\\
             &&-\,d\eta(Y,Z)g(X,W)+d\eta(Y,W)g(X,Z).
\end{eqnarray*}}
\end{lemma}
\proof Denoted with $\Phi$ the fundamental 2-form defined by $\Phi(X,Y):=g(X,\varphi Y)$, let us consider the derivation $R_{XY}$ of the tensor algebra $\mathfrak{T}(M)$, canonically induced from the $(1,1)$--tensor field $R(X,Y):=[\nabla_X,\nabla_Y]-\nabla_{[X,Y]}$. For any $X,Y,Z,W\in\Gamma(TM)$, we have:
{\setlength\arraycolsep{2pt}
\begin{eqnarray}
(R_{XY}\Phi)(Z,W)&=&R_{XY}(g(Z,\varphi W))-\Phi(R_{XY}Z,W)-\Phi(Z,R_{XY}W)\nonumber\\
                 &=&-\,g(R_{XY}Z,\varphi W)-g(R_{XY}\varphi Z,W)\label{eq151}\\
                 &=&-\,g(R(X,Y)Z,\varphi W)-g(R(X,Y)\varphi Z,W)\nonumber
\end{eqnarray}}Let us now compute again the term $(R_{XY}\Phi)(Z,W)$, using (\ref{eq303}). One has
{\setlength\arraycolsep{2pt}
\begin{eqnarray*}
(\nabla_X\nabla_Y\Phi)(Z,W)&=&X(\nabla_Y\Phi(Z,W))-\nabla_Y\Phi(\nabla_XZ,W)-\nabla_Y\Phi(Z,\nabla_XW)\\
                           &=&X(g(Z,(\nabla_Y\varphi)(W)))-g(\nabla_XZ,(\nabla_Y\varphi)(W))+g((\nabla_Y\varphi)(Z),\nabla_XW)\\
                           &=&r\varepsilon\left(X(\eta(Z)g(Y,W))-X(\eta(W)g(Z,Y))-\eta(\nabla_XZ)g(Y,W)\right.\\
                            &&+\left.\eta(W)g(\nabla_XZ,Y)+\eta(\nabla_XW)g(Y,Z)-\eta(Z)g(Y,\nabla_XW)\right)
\end{eqnarray*}}Switching $X$ and $Y$, we have
{\setlength\arraycolsep{2pt}
\begin{eqnarray*}
(\nabla_Y\nabla_X\Phi)(Z,W)&=&r\varepsilon\left(Y(\eta(Z)g(X,W))-Y(\eta(W)g(Z,X))-\eta(\nabla_YZ)g(X,W)\right.\\
                            &&+\left.\eta(W)g(\nabla_YZ,X)+\eta(\nabla_YW)g(X,Z)-\eta(Z)g(X,\nabla_YW)\right).
\end{eqnarray*}}Finally
\[
(\nabla_{[X,Y]}\Phi)(Z,W)=g(Z,(\nabla_{[X,Y]}\varphi)(W))=r\varepsilon\left(\eta(Z)g([X,Y],W)-\eta(W)g(Z,[X,Y])\right).                         
\]
It follows
{\setlength\arraycolsep{2pt}
\begin{eqnarray*}
(R_{XY}\Phi)(Z,W)&=&r\varepsilon\left((\nabla_X\eta)(Z)g(Y,W)-(\nabla_X\eta)(W)g(Z,Y)\right.\\
                 &&-\left.(\nabla_Y\eta)(Z)g(X,W)+(\nabla_Y\eta)(W)g(Z,X)\right).
\end{eqnarray*}}Since $\nabla_X\xi=-\varepsilon\varphi(X)$ and $\Phi=d\eta$, then $(\nabla_X\eta)(Y)=d\eta(X,Y)$, and we have
{\setlength\arraycolsep{2pt}
\begin{eqnarray}
(R_{XY}\Phi)(Z,W)&=&r\varepsilon\left(d\eta(X,Z)g(Y,W)-d\eta(X,W)g(Z,Y)\right.\nonumber\\
                 &&-\left.d\eta(Y,Z)g(X,W)+d\eta(Y,W)g(Z,X)\right)\label{eq141}\\
                 &=&r\varepsilon P(X,Y,Z,W).\nonumber
\end{eqnarray}}Now, (\ref{eq151}) and (\ref{eq141}) imply (\ref{eq139}). \proofend

It is easy to prove the following Lemma.
\begin{lemma}\label{lemma010}
Let $M$ be a manifold endowed with an almost (para)contact metric structure $(\varphi,\xi,\eta,g)$. Then one has, for any $X_1,X_2,X_3,X_4\in\Gamma(TM)$:
\begin{enumerate}
 \item[\emph{(i)}]$P(X_1,X_2,X_3,X_4)=-P(X_2,X_1,X_3,X_4)$;
 \item[\emph{(ii)}]$P(X_1,X_2,X_3,X_4)=-P(X_1,X_2,X_4,X_3)$;
 \item[\emph{(iii)}]$P(X_1,X_2,X_3,X_4)=-P(X_3,X_4,X_1,X_2)$;
 \item[\emph{(iv)}]$P(X_1,X_2,X_3,X_4)=P(X_4,X_3,X_2,X_1)$;
\end{enumerate}
\end{lemma}

\begin{proposition}
Let $M^{2n+1}$ be a manifold with an $[r]$-Sasakian structure $(\varphi,\xi,\eta,g)$. Then
\begin{equation}\label{eq142}
\rho(X,\xi)=2nr\eta(X),
\end{equation}
for any $X\in\Gamma(TM)$.
\end{proposition}
\proof We choose a local orthonormal frame $(E_i)_{1\leqslant i\leqslant 2n+1}$ on $M$. Putting $\alpha_i:=g(E_i,E_i)$, using (\ref{Eq303}), (\ref{eq139}) and the definition of $P$, since $I=-r\varphi^2+\eta\otimes\xi$ and $d\eta(X,Y)=\Phi(X,Y)=g(X,\varphi Y)$, one has, for any $X\in\Gamma(TM)$:
{\setlength\arraycolsep{2pt}
\begin{eqnarray*}
\rho(X,\xi)&=&\sum_{i=1}^{2n+1}\alpha_iR(X,E_i,\xi,E_i)=-r\sum_{i=1}^{2n+1}\alpha_iR(X,E_i,\xi,\varphi^2E_i)\\
           &=&-r\left(\sum_{i=1}^{2n+1}\alpha_ig(R(X,E_i)\varphi(\xi),\varphi E_i)+\varepsilon r\sum_{i=1}^{2n+1}\alpha_iP(X,E_i,\xi,\varphi E_i)\right)\\
           &=&-\varepsilon\sum_{i=1}^{2n+1}\alpha_i\left(g(\varphi X,\varphi E_i)g(\xi,E_i)-g(\varphi E_i,\varphi E_i)g(X,\xi)\right)\\
           &=&-r\varepsilon\sum_{i=1}^{2n+1}\alpha_i\left(g(X,E_i)g(\xi,E_i)-g(E_i,E_i)g(X,\xi)\right)\\
           &=&-r\varepsilon\left\{g(X,\xi)-\sum_{i=1}^{2n+1}\alpha_i^2g(X,\xi)\right\}=r2n\eta(X).
\end{eqnarray*}}\proofend

\section{Mixed 3-Sasakian structures and Ricci curvature}

Let us state now the main result, and then examine few consequences.

\begin{theorem}\label{Propo001}
Any mixed 3-Sasakian manifold $(M^{4n+3},\varphi_a,\xi_a,\eta^a,g)$ is Einstein. More precisely, for any $X,Y\in\Gamma(TM)$, one has
\[
\rho(X,Y)=-\sigma(4n+2)g(X,Y),
\]
where $\sigma=\pm1$, according as the 3-structure is positive or negative.
\end{theorem}
\proof Let us put, for any $X,Y\in\Gamma(TM)$:
\begin{equation}\label{eq155}
Q(X,Y):=\rho(X,\varphi_3Y)-\rho(Y,\varphi_3X)+2\sigma(4n+1)g(X,\varphi_3Y).
\end{equation}
We are going to prove that
\begin{equation}\label{Eq005}
Q(X,Y)=\sum_{i=1}^{4n+3}\varepsilon_ig(R(X,Y)e_i,\varphi_3(e_i)),
\end{equation}
where $(e_i)_{1\leqslant i\leqslant 4n+3}$ is an arbitrary orthonormal local frame on $M$, and $\varepsilon_i:=g(e_i,e_i)$. Since the structure $(\varphi_3,\xi_3,\eta^3,g)$ is indefinite Sasakian, from (\ref{eq139}), with $r=1$ and $\varepsilon=g(\xi_3,\xi_3)=\mp1=-\sigma$, according as the 3-structure is positive or negative, we have
\begin{equation}\label{eq152}
g(R(X,Y)Z,\varphi_3W)=-g(R(X,Y)\varphi_3Z,W)+\sigma P_3(X,Y,Z,W)
\end{equation}
for any $X,Y,Z,W\in\Gamma(TM)$.

Using Bianchi's First Identity, (\ref{eq152}), and Lemma \ref{lemma010}, the right hand side of (\ref{Eq005}) becomes
{\setlength\arraycolsep{2pt}
\begin{eqnarray}
\sum_{i=1}^{4n+3}\varepsilon_ig(R(X,Y)e_i,\varphi_3(e_i))&=&-\sum_{i=1}^{4n+3}\varepsilon_i\left\{g(R(Y,e_i)X,\varphi_3(e_i))+g(R(e_i,X)Y,\varphi_3(e_i))\right\}\nonumber \\
      &=&\sum_{i=1}^{4n+3}\varepsilon_i\left\{g(R(Y,e_i)\varphi_3X,e_i)-\sigma P_3(Y,e_i,X,e_i)\right.\label{Eq006} \\
       &&+\left.g(R(e_i,X)\varphi_3Y,e_i)-\sigma P_3(e_i,X,Y,e_i)\right\}\nonumber \\
      &=&-\,\rho(Y,\varphi_3X)+\rho(X,\varphi_3Y)-2\sigma\sum_{i=1}^{4n+3}\varepsilon_iP_3(Y,e_i,X,e_i).\nonumber
\end{eqnarray}}Computing the last term, by the definition of $P_3$, one has
{\setlength\arraycolsep{2pt}
\begin{eqnarray}
\sum_{i=1}^{4n+3}\varepsilon_iP_3(Y,e_i,X,e_i)&=&\sum_{i=1}^{4n+3}\varepsilon_i\left\{d\eta^3(Y,X)g(e_i,e_i)-d\eta^3(Y,e_i)g(e_i,X)\right.\nonumber\\
                                                 &&-\left.d\eta^3(e_i,X)g(Y,e_i)-d\eta^3(e_i,e_i)g(Y,X)\right\}\nonumber\\
                                              &=&(4n+3)g(\varphi_3X,Y)+g(X,\varphi_3Y)-g(\varphi_3X,Y)\label{eq154}\\
                                              &=&(4n+3)g(\varphi_3X,Y)-2g(\varphi_3X,Y)\nonumber\\
                                              &=&-(4n+1)g(X,\varphi_3Y).\nonumber
\end{eqnarray}}From (\ref{Eq006}) and (\ref{eq154}), we obtain (\ref{eq155}).

Now, we choose a local orthonormal frame adapted to the 3-structure of type
\[
(E_i,\varphi_1E_i,\varphi_2E_i,\varphi_3E_i,\xi_1,\xi_2,\xi_3)_{1\leqslant i\leqslant n}.
\]
We put, for any $i\in\{1,\ldots,n\}$: $e_i:=E_i$, $e_{n+i}:=\varphi_1E_i$, $e_{2n+i}:=\varphi_2E_i$ and $e_{3n+i}:=\varphi_3E_i$, and
\[
\alpha_i:=g(E_i,E_i)=-g(\varphi_1E_i,\varphi_1E_i)=-g(\varphi_2E_i,\varphi_2E_i)=g(\varphi_3E_i,\varphi_3E_i);
\]
for any $a\in\{1,2,3\}$, we put also: $e_{4n+a}:=\xi_a$, and $\alpha_{4n+a}:=g(\xi_a,\xi_a)=\varepsilon_a$. Then we get
{\setlength\arraycolsep{2pt}
\begin{eqnarray*}
Q(X,Y)&=&\sum_{i=1}^{n}\alpha_i\{g(R(X,Y)E_i,\varphi_3E_i)-g(R(X,Y)\varphi_1E_i,\varphi_3\varphi_1E_i)\\
                                  &&-\,g(R(X,Y)\varphi_2E_i,\varphi_3\varphi_2E_i)+g(R(X,Y)\varphi_3E_i,\varphi^2_3E_i)\}\\
                                  &&+\,\varepsilon_1g(R(X,Y)\xi_1,\varphi_3\xi_1)+\varepsilon_2g(R(X,Y)\xi_2,\varphi_3\xi_2)\\
                                  &=&\sum_{i=1}^{n}\alpha_i\{g(R(X,Y)E_i,\varphi_3E_i)+g(R(X,Y)\varphi_1E_i,\varphi_1\varphi_3E_i)\\
                                  &&+\,g(R(X,Y)\varphi_2E_i,\varphi_2\varphi_3E_i)+g(R(X,Y)E_i,\varphi_3E_i)\}\\
                                  &&+\,\varepsilon_1g(R(X,Y)\xi_1,\varphi_1\xi_3)+\varepsilon_2g(R(X,Y)\xi_2,\varphi_2\xi_3).
\end{eqnarray*}}Since the structures $(\varphi_1,\xi_1,\eta^1,g)$ and $(\varphi_2,\xi_2,\eta^2,g)$ are both paraSasakian, using (\ref{eq139}) with $r=-1$, one has
{\setlength\arraycolsep{2pt}
\begin{eqnarray*}
Q(X,Y)&=&\sum_{i=1}^{n}\alpha_i\left\{g(R(X,Y)E_i,\varphi_3E_i)-g(R(X,Y)\varphi_1^2E_i,\varphi_3E_i)+\varepsilon_1P_1(X,Y,\varphi_1E_i,\varphi_3E_i)\right.\\
         &&-\left.g(R(X,Y)\varphi_2^2E_i,\varphi_3E_i)+\varepsilon_2P_2(X,Y,\varphi_2E_i,\varphi_3E_i)+g(R(X,Y)E_i,\varphi_3E_i)\right\}\\
         &&+\,P_1(X,Y,\xi_1,\xi_3)+P_2(X,Y,\xi_2,\xi_3)\\
       &=&\sum_{i=1}^{n}\alpha_i\{\varepsilon_1P_1(X,Y,\varphi_1E_i,\varphi_3E_i)+\varepsilon_2P_2(X,Y,\varphi_2E_i,\varphi_3E_i)\}\\
         &&+\,P_1(X,Y,\xi_1,\xi_3)+P_2(X,Y,\xi_2,\xi_3).
\end{eqnarray*}}Recalling the definition of the tensor field $P$, since $d\eta^1=\Phi_1$, $d\eta^2=\Phi_2$, $\varepsilon_1=\varepsilon_2=\sigma=-\varepsilon_3$ and $\sigma\varepsilon_1=\sigma\varepsilon_2=1$, using (\ref{eq90}), (\ref{eq92}) and (\ref{eq93}), one has
{\setlength\arraycolsep{2pt}
\begin{eqnarray}
Q(X,Y)&=&-2\sigma\bigg\{\sum_{i=1}^{n}\alpha_i((g(X,E_i)g(\varphi_3Y,E_i)-g(X,\varphi_2E_i)g(\varphi_3Y,\varphi_2E_i)\nonumber\\
      &&+\,g(\varphi_3Y,\varphi_3E_i)g(X,\varphi_3E_i)-g(\varphi_3Y,\varphi_1E_i)g(X,\varphi_1E_i))\label{eq156}\\
      &&+\,\varepsilon_1g(X,\xi_1)g(\varphi_3Y,\xi_1)+\varepsilon_2g(X,\xi_2)g(\varphi_3Y,\xi_2)\bigg\}\nonumber\\
      &=&-2\sigma g(X,\varphi_3Y).\nonumber
\end{eqnarray}}From (\ref{eq155}) and (\ref{eq156}), it follows
\begin{equation}\label{eq157}
\rho(X,\varphi_3Y)-\rho(\varphi_3X,Y)=-2\sigma(4n+2)g(X,\varphi_3Y).
\end{equation}
Since the structure $(\varphi_3,\xi_3,\eta^3,g)$ is indefinite Sasakian, if $X,Y$ are orthogonal to $\xi_3$, one has (cf.\ \cite{Bla} for the Riemannian case) $\rho(X,\varphi_3Y)=-\rho(\varphi_3X,Y)$, hence, from (\ref{eq157}), for any $X,Y$ orthogonal to $\xi_3$, it follows that $\rho(X,\varphi_3Y)=-\sigma(4n+2)g(X,\varphi_3Y)$. Replacing $Y$ with $\varphi_3Y$, since $Y$ is orthogonal to $\xi_3$, one has
\begin{equation}\label{eq158}
\rho(X,Y)=-\sigma(4n+2)g(X,Y), \qquad X,Y\bot\,\xi_3.
\end{equation}
Using (\ref{eq142}), we have
\begin{equation}\label{eq159}
\rho(X,\xi_3)=-\sigma(4n+2)g(X,\xi_3), \qquad X\in\Gamma(TM),
\end{equation}
hence, in particular putting $X=\xi_3$:
\begin{equation}\label{eq160}
\rho(\xi_3,\xi_3)=-\sigma(4n+2)g(\xi_3,\xi_3).
\end{equation}
Finally, if $X,Y\in\Gamma(TM)$, putting $X=X_0+\lambda\xi_3$ and $Y=Y_0+\mu\xi_3$, with $X_0,Y_0$ orthogonal to $\xi_3$, and $\lambda,\mu\in\mathfrak{F}(M)$, using (\ref{eq158}), (\ref{eq159}) and (\ref{eq160}), one gets $\rho(X,Y)=-\sigma(4n+2)g(X,Y)$ for any $X,Y\in\Gamma(TM)$, concluding the proof. \proofend

As an obvious consequence of the above result, we have.
\begin{proposition}\label{Propo002}
Any mixed 3-Sasakian manifold $(M^{4n+3},\varphi_a,\xi_a,\eta^a,g)$ has constant scalar curvature
\[
\mathrm{Sc}=-\sigma(4n+2)(4n+3),
\]
therefore negative or positive, according as the 3-structure is positive or negative.
\end{proposition}

\begin{proposition}\label{Propo003}
Let $(M^{4n+3},\varphi_a,\xi_a,\eta^a,g)$ be a mixed 3-Sasakian manifold. Then, $M$ has (pointwise) constant sectional curvature $k$ if and only if $k=\mp1$, according as the 3-structure is positive or negative, respectively.
\end{proposition}
\proof Since the 3-structure \mmts{ }is mixed 3-Sasakian, then (\ref{eq139}) holds for any $a\in\{1,2,3\}$. Using the constant $\sigma=\pm1$, according as the 3-structure is positive or negative, and recalling that $\tau_a\varepsilon_a=-\sigma$, we have, for any $a\in\{1,2,3\}$ and any $X,Y,Z,W\in\Gamma(TM)$, $g(R(X,Y)Z,\varphi_aW)+g(R(X,Y)\varphi_aZ,W)=\sigma P_a(X,Y,Z,W)$. Supposing that $M$ has pointwise constant sectional curvature $k\in\mathfrak{F}(M)$, i.e. $R(X,Y)Z=k\{g(Y,Z)X-g(X,Z)Y\}$, we have
{\setlength\arraycolsep{2pt}
\begin{eqnarray*}
\sigma P_a(X,Y,Z,W)&=&g(R(X,Y)Z,\varphi_aW)+g(R(X,Y)\varphi_aZ,W)\\
                   &=&k\{g(Y,Z)g(X,\varphi_aW)-g(X,Z)g(Y,\varphi_aW)\\
                     &&+\,g(Y,\varphi_aZ)g(X,W)-g(X,\varphi_aZ)g(Y,W)\}\\
                   &=&k\{d\eta^a(X,W)g(Y,Z)-d\eta^a(Y,W)g(X,Z)\\
                     &&+\,d\eta^a(Y,Z)g(X,W)-d\eta^a(X,Z)g(Y,W)\}\\
                   &=&-\,kP_a(X,Y,Z,W),
\end{eqnarray*}}hence, for any $a\in\{1,2,3\}$ and any $X,Y,Z,W\in\Gamma(TM)$, we get $(k+\sigma)P_a(X,Y,Z,W)=0$, and it follows that $k=-\sigma=\mp1$, according as the 3-structure is positive or negative. Namely, choosing a vector field $Y$ orthogonal to $\xi_1,\xi_2,\xi_3$ and such that $g(Y,Y)\neq0$, by the definition of $P_a$, given in Lemma \ref{lemma139}, we get $P_a(\xi_a,Y,\xi_a,\varphi_aY)=-\varepsilon_ag(Y,Y)\neq0$. \proofend

\section{Mixed metric 3-contact and mixed 3-Sasakian structures}

In this section we shall be concerned with some properties of particular classes of mixed metric 3-structures, namely the class of the mixed metric 3-contact structures, which reflect analogous properties of classical metric 3-structures; for more details see \cite{Bla}.

\begin{definition}
Let $M$ be a manifold with a mixed metric 3-structure \mmts. This structure is said to be a \emph{mixed metric 3-contact structure} if $d\eta^a=\Phi_a$, for each $a\in\{1,2,3\}$, where $\Phi_a$ is the fundamental 2-form defined by $\Phi_a(X,Y):=g(X,\varphi_aY)$. Then $(M,\varphi_a,\xi_a,\eta^a,g)$ will be said \emph{mixed metric 3-contact manifold}.
\end{definition}

Our intent here is to prove that any mixed metric 3-contact manifold is in fact a mixed 3-Sasakian manifold.  

Let $M$ be a manifold with a mixed metric 3-structure \mmts. Setting $\tilde{M}=M\times\R$, and denoting with $t$ the coordinate on $\R$, define three $(1,1)$--type tensor fields $J_a$, ${a=1,2,3}$, by putting, for any $\tilde{X}=\left(X,f\frac{d}{dt}\right)\in\Gamma(T\tilde{M})$, with $X\in\Gamma(TM)$ and $f\in\mathfrak{F}(\tilde{M})$:
\[
J_a(\tilde{X})=J_a\left(X,f\frac{d}{dt}\right):=\left(\varphi_aX-\tau_af\xi_a,\eta^a(X)\frac{d}{dt}\right),
\]
where $\tau_1=\tau_2=-1=-\tau_3$. Furthermore, define the $(0,2)$--type tensor field $G$, by putting, for any $\tilde{X}=\left(X,f\frac{d}{dt}\right)$ and $\tilde{Y}=\left(Y,h\frac{d}{dt}\right)$ in $\Gamma(T\tilde{M})$, with $X,Y\in\Gamma(TM)$ and $f,h\in\mathfrak{F}(\tilde{M})$:
\[
G(\tilde{X},\tilde{Y}):=g(X,Y)-\sigma fh,
\]
where $\sigma=\pm1$, according as the 3-structure is positive or negative, respectively. 

\begin{proposition}
$(\tilde{M},J_a,G)_{a=1,2,3}$ is an almost hyper paraHermitian manifold.
\end{proposition}
\proof Let be $a\in\{1,2,3\}$, and $\tilde{X}\in\Gamma(T\tilde{M})$, with $\tilde{X}=\left(X,f\frac{d}{dt}\right)$. Since by definition $\varphi_a^2=-\tau_a(I-\eta^a\otimes\xi_a)$, we have
{\setlength\arraycolsep{2pt}
\begin{eqnarray*}
(J_a)^2(\tilde{X}) =\left((\varphi_a)^2X-\tau_a\eta^a(X)\xi_a,-\tau_af\frac{d}{dt}\right)
                =-\tau_a\tilde{X},
\end{eqnarray*}}hence $(J_a)^2=-\tau_aI$. Let now $(a,b,c)$ be a cyclic permutation of $(1,2,3)$. Using (\ref{eq90}), (\ref{eq91}) and (\ref{eq92}), one has, for any $\tilde{X}\in\Gamma(T\tilde{M})$, with $\tilde{X}=\left(X,f\frac{d}{dt}\right)$:
{\setlength\arraycolsep{2pt}
\begin{eqnarray*}
J_aJ_b(\tilde{X})&=&\left(\varphi_a\varphi_bX-\tau_bf\varphi_a\xi_b-\tau_a\eta^b(X)\xi_a,(\eta^a(\varphi_bX)-\tau_bf\eta^a\xi_b)\frac{d}{dt}\right)\\
                 &=&\left(\tau_c\varphi_cX-f\xi_c,\tau_c\eta^c(X)\frac{d}{dt}\right)=\tau_cJ_c(\tilde{X})
\end{eqnarray*}}hence $J_aJ_b=\tau_cJ_c$. Analogously, one has $J_bJ_a=-\tau_cJ_c$, and this proves that $(J_a)_{a=1,2,3}$ is an almost hyper paracomplex structure on $\tilde{M}$. Let now be $a\in\{1,2,3\}$, $\tilde{X}=\left(X,f\frac{d}{dt}\right)$ and $\tilde{Y}=\left(Y,h\frac{d}{dt}\right)$. Since, by (\ref{eq93}), $g(\varphi_aX,Y)=-g(X,\varphi_aY)$, using the identity $\tau_a\varepsilon_a=-\sigma$, by standard calculations we have $G(\tilde{X},J_a\tilde{Y})=-G(J_a(\tilde{X}),\tilde{Y})$ and, by Definition \ref{def005} it follows that $(\tilde{M},J_a,G)$, $a\in\{1,2,3\}$, is an almost hyper paraHermitian manifold. \endproof

\begin{remark}
It is clear that, for any $a\in\{1,2,3\}$ the tensor field $J_a$ constructed on $\tilde{M}$ is an almost product structure, for $a=1,2$, and an almost complex structure, for $a=3$. The three structures \mmts{ }are normal if and only if the manifold $(\tilde{M},J_a,G)$, $a\in\{1,2,3\}$, is hyper paraHermitian.
\end{remark}

Thus, we may state:

\begin{proposition}
Let $M$ be a manifold endowed with a mixed 3-structure \mts. Then, the structures are normal if and only if at least two of them are normal.
\end{proposition}

We shall see in a moment that the manifold $(\tilde{M},J_a,G)$, $a\in\{1,2,3\}$, is indeed hyper paraHermitian, if the 3-structure is a mixed metric 3-contact structure. To this purpose, let us prove the following preliminary results.

\begin{lemma}\label{lemma006}
Let $M$ be a manifold endowed with a mixed metric 3-contact structure. Denoting, for any $a\in\{1,2,3\}$, with $\Omega_a$ the fundamental 2-form associated with the structure $(J_a,G)$, defined by $\Omega_a(\tilde{X},\tilde{Y}):=G(\tilde{X},J_a\tilde{Y})$, we have
\[
d\Omega_a=2\sigma dt\wedge\Omega_a,
\]
for any $a\in\{1,2,3\}$, where $\sigma=\pm1$ according as the 3-structure is positive or negative.
\end{lemma}
\proof Fixed $a\in\{1,2,3\}$, let us compute $d\Omega_a$, using the formula
\begin{equation}\label{eq132}
3d\Omega_a(\tilde{X},\tilde{Y},\tilde{Z})=\stackunder{\scriptscriptstyle{(\tilde{X},\tilde{Y},\tilde{Z})}}{\mathfrak{S}}\{\tilde{X}(\Omega_a(\tilde{Y},\tilde{Z}))-\Omega_a([\tilde{X},\tilde{Y}],\tilde{Z})\},
\end{equation}
for any $\tilde{X},\tilde{Y},\tilde{Z}\in\Gamma(T\tilde{M})$. To this purpose, putting $\tilde{X}=\left(X,f\frac{d}{dt}\right)$, $\tilde{Y}=\left(Y,h\frac{d}{dt}\right)$ and $\tilde{Z}=\left(Z,k\frac{d}{dt}\right)$, using $\tau_a\varepsilon_a=-\sigma$ we have
\begin{equation}\label{eq131}
\Omega_a(\tilde{Y},\tilde{Z})=\Phi_a(Y,Z)+\sigma(k\eta^a(Y)-h\eta^a(Z)).
\end{equation}
Furthermore $[\tilde{X},\tilde{Y}]=\bigg([X,Y],\left(X(h)-Y(f)+f\frac{dh}{dt}-h\frac{df}{dt}\right)\frac{d}{dt}\bigg)$ and
\[
\Omega_a([\tilde{X},\tilde{Y}],\tilde{Z})=\Phi_a([X,Y],Z)+\sigma\bigg\{k\eta^a[X,Y]-\left(X(h)-Y(f)+f\frac{dh}{dt}-h\frac{df}{dt}\right)\eta^a(Z)\bigg\}.
\]
Finally, from (\ref{eq131})
{\setlength\arraycolsep{2pt}
\begin{eqnarray*}
\tilde{X}(\Omega_a(\tilde{Y},\tilde{Z}))&=&X(\Phi_a(Y,Z)+\sigma(k\eta^a(Y)-h\eta^a(Z)))\\
                                        &&+\,f\frac{d}{dt}(\Phi_a(Y,Z)+\sigma(k\eta^a(Y)-h\eta^a(Z)))\\
                                        &=&X(\Phi_a(Y,Z))+\sigma(X(k)\eta^a(Y)+kX(\eta^a(Y))-X(h)\eta^a(Z)-hX(\eta^a(Z)))\\
                                        &&+\,\sigma\left(f\frac{dk}{dt}\eta^a(Y)-f\frac{dh}{dt}\eta^a(Z)\right).
\end{eqnarray*}}From (\ref{eq132}), using the above identities, and $d\Phi_a=0$, one gets
\[
3d\Omega_a(\tilde{X},\tilde{Y},\tilde{Z})=2\sigma\left(\Phi_a(X,Y)k+\Phi_a(Y,Z)f+\Phi_a(Z,X)h\right).
\]
Finally, using (\ref{eq131}), it follows
{\setlength\arraycolsep{2pt}
\begin{eqnarray*}
3d\Omega_a(\tilde{X},\tilde{Y},\tilde{Z})&=&2\sigma\left(f\Omega_a(\tilde{Y},\tilde{Z})-\sigma(fk\eta^a(Y)-fh\eta^a(Z))\right.\\
                                            &&+\left.h\Omega_a(\tilde{Z},\tilde{X})-\sigma(hf\eta^a(Z)-hk\eta^a(X))\right.\\
                                            &&+\left.k\Omega_a(\tilde{X},\tilde{Y})-\sigma(kh\eta^a(X)-kf\eta^a(Y))\right)\\
                                            &=&2\sigma(f\Omega_a(\tilde{Y},\tilde{Z})+h\Omega_a(\tilde{Z},\tilde{X})+k\Omega_a(\tilde{X},\tilde{Y}))\\
                                            &=&6\sigma(dt\wedge\Omega_a)(\tilde{X},\tilde{Y},\tilde{Z}),
\end{eqnarray*}}hence $d\Omega_a=2\sigma dt\wedge\Omega_a$. \endproof

\begin{lemma}\label{lemma007}
Let $(M,J_a,g)$, $a\in\{1,2,3\}$, be an almost hyper paraHermitian manifold, such that, denoting with $\Omega_a$ the fundamental 2-form associated with $J_a$, there exists a 1-form $\omega$ satisfying $d\Omega_a=k\omega\wedge\Omega_a$, for any $a\in\{1,2,3\}$, with $k\in\mathfrak{F}(M)$. Then, each structure $J_a$ is integrable, and the manifold is hyper paraHermitian.
\end{lemma}
\proof Let us prove that $N_1=0$. It is well known that
\begin{equation*}
N_1(X,Y)=(\nabla_{J_1X}J_1)(Y)-(\nabla_{J_1Y}J_1)(X)-J_1(\nabla_XJ_1)(Y)+J_1(\nabla_YJ_1)(X),
\end{equation*}
hence, using (i) and (ii) of Definition \ref{def005}, we get
\begin{equation}\label{eq138}
J_2N_1(X,Y)=-J_2(\nabla_{J_1Y}J_1)(X)-J_3(\nabla_YJ_1)(X)+J_2(\nabla_{J_1X}J_1)(Y)+J_3(\nabla_XJ_1)(Y).
\end{equation}
Then, for any $Z\in\Gamma(TM)$, using (iii) of Definition \ref{def005}, with standard calculations, one has
{\setlength\arraycolsep{2pt}
\begin{eqnarray*}
g(-J_2(\nabla_{J_1Y}J_1)(X),Z)&=&-g(J_2\nabla_{J_1Y}(J_1X),Z)-g(J_3\nabla_{J_1Y}X,Z)\\
                              &=&-g(X,(\nabla_{J_1Y}J_3)(Z))-g(J_1X,(\nabla_{J_1Y}J_2)(Z))\\
                              &=&(\nabla_{J_1Y}\Omega_3)(Z,X)+(\nabla_{J_1Y}\Omega_2)(Z,J_1X).
\end{eqnarray*}}Switching $X$ and $Y$ one has $g(J_2(\nabla_{J_1X}J_1)(Y),Z)=(\nabla_{J_1X}\Omega_3)(Y,Z)+(\nabla_{J_1X}\Omega_2)(J_1Y,Z)$. Analogously, one obtains $g(-J_3(\nabla_YJ_1)(X),Z)=(\nabla_Y\Omega_2)(Z,X)+(\nabla_Y\Omega_3)(Z,J_1X)$ and switching $X$ and $Y$ one gets $g(J_3(\nabla_XJ_1)(Y),Z)=(\nabla_X\Omega_2)(Y,Z)+(\nabla_X\Omega_3)(J_1Y,Z)$. Since $3d\Omega(X,Y,Z)=\stackunder{(X,Y,Z)}{\mathfrak{S}}(\nabla_X\Omega)(Y,Z)$, from (\ref{eq138}) we have
{\setlength\arraycolsep{2pt}
\begin{eqnarray*}
g(J_2N_1(X,Y),Z)&=&3d\Omega_2(X,Y,Z)+3d\Omega_3(X,J_1Y,Z)\\
                &&+\,3d\Omega_3(J_1X,Y,Z)+3d\Omega_2(J_1X,J_1Y,Z).
\end{eqnarray*}}Being $d\Omega_a=k\omega\wedge\Omega_a$, we get
{\setlength\arraycolsep{0.5pt}
\begin{eqnarray*}
g(J_2N_1(X,&&Y),Z)=k\left\{\omega(X)\Omega_2(Y,Z)+\omega(Y)\Omega_2(Z,X)+\omega(Z)\Omega_2(X,Y)\right.\\
                  &&\qquad\qquad+\left.\omega(X)\Omega_3(J_1Y,Z)+\omega(J_1Y)\Omega_3(Z,X)+\omega(Z)\Omega_3(X,J_1Y)\right.\\
                  &&\qquad\qquad+\left.\omega(J_1X)\Omega_3(Y,Z)+\omega(Y)\Omega_3(Z,J_1X)+\omega(Z)\Omega_3(J_1X,Y)\right.\\
                  &&\qquad\qquad+\left.\omega(J_1X)\Omega_2(J_1Y,Z)+\omega(J_1Y)\Omega_2(Z,J_1X)+\omega(Z)\Omega_2(J_1X,J_1Y)\right\}.
\end{eqnarray*}}It is easy to check that $\Omega_3(J_1Y,Z)=-\Omega_2(Y,Z)$, $\Omega_3(Y,J_1Z)=-\Omega_2(Y,Z)$, $\Omega_2(Z,J_1X)=-\Omega_3(Z,X)$, $\Omega_2(J_1Z,X)=-\Omega_3(Z,X)$ and $\Omega_2(J_1X,J_1Y)=\Omega_2(X,Y)$. Therefore, we obtain $g(J_2N_1(X,Y),Z)=0$, hence $N_1=0$. In an analogous way, one proves that $N_2=0$ and $N_3=0$. \endproof

As an obvious consequence of Lemmas \ref{lemma006} and \ref{lemma007}, one obtains the following result.

\begin{theorem}
Any mixed metric 3-contact structure on a manifold is mixed 3-Sasakian.
\end{theorem}

Thus, Theorems \ref{Propo001}, \ref{Propo002} and \ref{Propo003} may be formulated for mixed metric 3-contact manifolds.

\vspace{10pt}
\noindent \textbf{Acknowledgements.} The authors are grateful to Prof.\ Dr.\ S. Ianu\c{s} for many discussions about the topic of this paper, occurred during his stay at the University of Bari and during the first author's visit at the University of Bucharest.

\vspace{12pt}
\noindent\small{Authors' address}
\vspace{2pt}\\
  \small{Department of Mathematics, University of Bari} \\
  \small{Via E. Orabona 4,}\\
  \small{I-70125 Bari (Italy)}
\vspace{2pt}\\
  \small{\texttt{caldarella@dm.uniba.it}}\\
  \small{\texttt{pastore@dm.uniba.it}}

\end{document}